\documentclass[12pt,a4paper]{article}
\usepackage{mathrsfs}

\usepackage{amssymb}

\usepackage{amsmath}

\newtheorem{theorem}{Theorem}[section]
\newtheorem{definition}[theorem]{Definition}
\newtheorem{lemma}[theorem]{Lemma}
\newtheorem{corollary}[theorem]{Corollary}

\begin{document}
\title{Gr\"{o}bner-Shirshov bases for categories\footnote{Supported by the
NNSF of China (Nos.10771077, 10911120389).}}
\author{
L. A. Bokut\footnote {Supported by RFBR 09-01-00157, LSS--3669.2010.1
and SB RAS Integration grant No. 2009.97 (Russia).} \\
{\small \ School of Mathematical Sciences, South China Normal
University}\\
{\small Guangzhou 510631, P. R. China}\\
{\small Sobolev Institute of Mathematics, Russian Academy of
Sciences}\\
{\small Siberian Branch, Novosibirsk 630090, Russia}\\
{\small  bokut@math.nsc.ru}\\
\\
Yuqun Chen\footnote {Corresponding author.} \ and Yu Li
\\
{\small \ School of Mathematical Sciences, South China Normal
University}\\
{\small Guangzhou 510631, P. R. China}\\
{\small  yqchen@scnu.edu.cn}\\
{\small liyu820615@126.com}}

\date{}

\maketitle \noindent\textbf{Abstract:} In this paper we establish
Composition-Diamond lemma  for small categories. We give
Gr\"{o}bner-Shirshov bases for simplicial category and cyclic
category.



\section{Introduction}

This paper devotes to Gr\"{o}bner-Shirshov bases for small
categories (all categories below are supposed to be small) presented
by a graph (=quiver) and defining relations (see, Maclane \cite{Sa
M}). As important examples, we use the simpicial and the cyclic
categories (see, for example, Maclane \cite{Sa}, Gelfand, Manin
\cite{G}). In an above presentation, a category is viewed as a
``monoid with several objects". A free category $C(X)$, generated by
a graph $X$, is just ``free partial monoid of partial words"
$u=x_{i_1}\dots x_{i_n},\ n\geq 0,\ x_i\in X$ and all product
defined in $C(X)$. A relation is an expression $u=v,\ u,v\in C(X)$,
where sources and targets of $u,v$ are coincident respectively. The
same as for semigroups, we may use two equivalent languages:
Gr\"{o}bner-Shirshov bases language and rewriting systems language.
Since we are using the former,  we need a Composition-Diamond lemma
(CD-lemma for short) for a free associative partial  algebra $kC(X)$
over a field $k$, where $kC(X)$ is just a linear combination of
uniform (with the same sources and targets) partial words. Then it
is a routing matter to establish CD-lemma for $kC(X)$. It is a free
category (``semigroup") partial algebra of a free category. Remark
that in the literature one is usually used a language of rewriting
system, see, for example, Malbos \cite{Malbos}. Let us stress that
the partial associative algebras presented by graphs and defining
relations are closely related to the well known quotients of the
path algebras from representation theory of finitely dimensional
algebras, see, for example, Assem, Simson, Skowro\'{n}ski
\cite{ASS}. In this respect Gr\"{o}bner-Shirshov bases for
categories are closely related to non-commutative Gr\"{o}bner bases
for quotients of path algebras, see Farkas, Feustel, Green
\cite{ffg93}. Rewriting system language for non-commutative
Gr\"{o}bner bases of quotients of path algebras was used by
Kobayashi \cite{kobayashi2006}. Main new results of this paper are
Gr\"{o}bner-Shirshov bases for the simplicial and cyclic categories.

 All algebras are assumed to be over a field.

\section{A short survey on Gr\"{o}bner-Shirshov bases}

What is now called Gr\"{o}bner and Gr\"{o}bner-Shirshov bases theory
was initiated by A. I. Shirshov (1921-1981) \cite{Sh62a,Sh62b}, 1962
for non-associative and Lie algebras, by H. Hironaka
\cite{Hi64,Hi64-}, 1964 for quotients of commutative infinite series
algebras (both formal and convergent), and by B. Buchberger
\cite{bu65,bu70}, 1965, 1970 for commutative algebras.

English translation of selected works of A. I. Shirshov, including
\cite{Sh62a,Sh62b}, is recently published \cite{Shir3}.

Remark that Shirshov's approach was a most universal as we
understand now since Lie algebra case becomes a model for many
classes of non-commutative and non-associative algebras (with
multiple operations), starting with associative algebras (see
below). Hironaka's papers on resolution of singularities of
algebraic varieties become famous very soon and Hironaka got Fields
Medal due to them few years latter. B. Buchberger's thesis
influenced very much many specialists in computer sciences, as well
as in commutative algebras and algebraic geometry, for huge
important applications of his bases, named him under his supervisor
W. Gr\"{o}bner (1898-1980).

Original Shirshov's approach for Lie algebras \cite{Sh62b}, 1962,
based on a notion of composition $[f,g]_w$ of two monic Lie
polynomials $f,g$ relative to associative word $w$, i.e., $f,g$ are
elements of a free Lie algebra $Lie(X)$ regarded as the subspace of
Lie polynomials of the free associative algebra $k\langle X\rangle$,
and $w\in X^*$, the free monoid generated by $X$. The definition of
Lie composition relies on a definition of associative composition
$(f,g)_w$ as (monic) associative polynomials (after worked out into
$f,g$ all Lie brackets $[x,y]=xy-yx)$ relative to
degree-lexicographical order on $X^*$. Namely, $(f,g)_w = fb-ag$,
where $w=acb,\ \bar f=ac,\ \bar g=cb,\ a,b,c\in X^*,\ c\neq 1$. Here
$\bar f$ means the leading (maximal) associative word of $f$. Then
$(f,g)_w$ belongs to associative ideal $Id(f,g)$ of $k\langle
X\rangle$ generated by $f,g$, and the leading word of $(f,g)_w$ is
less than $w$. Now we need to put some Lie brackets $[fb] - [ag]$ on
$fb-ag$ in such a way that the result would belong to Lie ideal
generated by $f,g$ (so we can not trouble bracketing into $f,g$) and
the leading associative monomial of $[fb]-[ag]$ must be less than
$w$. To overcome these obstacles Shirshov used his previous paper
\cite{sh58}, 1958 with a new linear basis of free Lie algebra
$Lie(X)$. As it happened the same linear basis of $Lie(X)$ was
discovered in the paper Chen, Fox, Lyndon \cite{CFL}, 1958. Now this
basis is called Lyndon-Shirshov basis, or, by a mistake, Lyndon
basis. It consists of non-associative Lyndon-Shirshov words (NLSW)
$[u]$ in $X$, that are in one-one correspondence with associative
Lyndon-Shirshov words (ALSW) $u$ in $X$. The latter is defined as by
a property $u=vw>wv$ for any $v,w\neq 1$. Shirshov \cite{sh58}, 1958
introduced and used the following properties of both associative and
non-associative Lyndon-Shirshov words:

1) For any ALSW u there is a unique bracketing $[u]$ such that $[u]$
is a NLSW.

There are two algorithms for bracketing an ALSW. He mostly used
``down-to-up algorithm"  to rewrite an ALSW u on a new alphabet
$X_u=\{x_i(x_\beta)^j,\ i>\beta,\ j\geq 0,\ x_{\beta} \mbox{ is the
minimal letter in } u\}$; the result $u'$ is again ALSW on $X_u$
with the lex-order $x_i>x_ix_{\beta}>((x_ix_\beta)x_{\beta})>\dots$.

It is Shirshov's rewriting or elimination algorithm from his  famous
paper
 Shirshov \cite{Shir53b}, 1953, on what is now called Shirshov-Witt theorem (any
subalgebra of a free Lie algebra is free). This rewriting was
rediscovered by Lazard \cite{Lazard60}, 1960 and now called as
Lazard elimination (it is better to call Lazard-Shirshov
elimination).

There is ``up-to-down algorithm" (see Shirshov \cite{sh58-}, 1958,
Chen, Fox, Lyndon \cite{CFL}, 1958): $[u]= [[v][w]]$, where $w$ is
the longest proper end of $u$ that is ALSW, in this case $v$ is also
an ALSW.

2) Leading associative word of NLSW $[u]$ is just $u$ (with the
coefficient 1).

3) Leading associative word of any Lie polynomial is an associative
Lyndon-Shirshov word.

4) A non-commutative polynomial $f$ is a Lie polynomial if and only
if
$$
f=f_0\rightarrow \cdots \rightarrow f_i\rightarrow
f_{i+1}\rightarrow \cdots \rightarrow f_n=0,
$$
where $f_i\rightarrow f_{i+1}=f_i-\alpha_i [u_i],\ \bar f_i= u_i$ is
an ALSW, $\alpha_i$ is the leading coefficient of $f_i,\
i=0,1,\dots$.

5) Any associative word $c\neq 1$ is the unique product of (not
strictly) increasing sequence of associative Lyndon-Shirshov words:
$c=c_1c_2\cdots c_n,\ c_1\leq \dots \leq c_n,\ c_i$ are ALSW's.

6) If $u=avb$, where $u,v$ are ALSW, $a,b\in X^*$, then there is a
relative bracketing $[u]_v=[a[v]b]$ of $u$ relative to $v$, such
that the leading associative word of $[u]_v$ is just $u$. Namely,
$[u]=[a[vc]d],\ cd=b,\ [u]_v=[a[[v[c_1]]\dots [c_n]]d],\
c=c_1c_2\cdots c_n$ as above.

7) If $ac$ and $cb$ are ASLW's and $c\neq 1$ then $acb$ is an ALSW
as well. If $a,b$ are ALSW's and $a>b$, then $ab$ is an ALSW as
well.

Property 5) was known to Chen, Fox, Lyndon \cite{CFL}, 1958 as well.
Lyndon \cite{Ly54}, 1954, was actually the first for definition of
associative ``Lyndon-Shirshov" words. To the best of our knowledge,
for many years, until PhD thesis by Viennot \cite{VG}, 1978, no one
mentioned the Lyndon's  discovery in 1954. On the other hand, there
were dozens of papers and some books on Lie algebras that mentioned
both associative and non-associative ``Lyndon-Shirshov" words as
Shirshov's regular words, see, for example, P. M. Cohn \cite{Conh},
1965, Bahturin \cite{Bah78}, 1978.

Now a Lie composition $[f,g]_w$ of monic Lie polynomials $f,g$
relative to a word $w=\bar fb=a\bar g =acb,\  c\neq 1$ is defined by
Shirshov \cite{Sh62b}, 1962, as follows
$$
[f,g]_w=[fb]_{\bar f} - [ag]_{\bar g},
$$
where $[fb]_{\bar f}$ means the result of substitution  f for $[\bar
f]$ into the relative bracketing  $[w]_{ac}$ of $w$ with respect to
$\bar f=ac$, the same for $[ag]_{\bar g}$.

According to the definition and properties above, any Lie
composition $[f,g]_w$ is an element of the Lie ideal generated by
$f,g$, and the leading associative word of the composition is less
than $w$.

The composition above is now called composition of intersection.
Shirshov avoided what is now called composition of inclusion
$$
[f,g]_w= f - [agb]_{\bar g},\ \ w=\bar f = a\bar gb,
$$
assuming that any system $S$ of Lie polynomials is reduced
(irreducible) in a sense that leading associative word  of any
polynomial from $S$ does not contain leading associative  words of
another polynomials from $S$. This assumption relies on his
algorithm of elimination of leading words for Lie polynomials below.

For associative polynomials the elimination algorithm is just
non-commutative version of the Euclidean elimination algorithm. For
Lie polynomial case Shirshov \cite{Sh62b}, 1962, defined the
elimination of a leading word as follows:

If $w=avb$, where $w, v$ are ALSW's, and $v$ is the leading word of
some monic Lie polynomial $f$, then the transformation $[w]\mapsto
[w] - [afb]_v$ is called an elimination of leading word of $f$ into
$[w]$. The result of Lie elimination is a Lie polynomial with a
leading associative word less than $w$.

Then Shirshov  \cite{Sh62b}, 1962, formulated an algorithm to add to
an initial reduced system of Lie polynomials $S$  a ``non-trivial"
composition $[f,g]_w$, where $f,g$ belong to $S$. Non-triviality  of
a Lie polynomial $h$ relative to $S$ means that  $h$ is not going to
zero using ``elimination of leading words of $S$". Actually, he
defines to add to $S$ not just a composition but rather the result
of elimination of leading words of $S$ into the composition in order
to have a reduced system as well.

Then Shirshov  proved the following

{\bf Composition Lemma.} Let $S$ be a reduced subset of $Lie(X)$. If
$f$ belongs to the Lie ideal generated by $S$, then the leading
associative word $\bar f$ contains, as a subword, some leading
associative word of a reduced multi-composition of elements of $S$.

He constantly used the following clear

{\bf Corollary.} The set of all irreducible NLSW's $[u]$ such that
$u$ does not contain any leading associative word of a reduced
multi-composition of elements of $S$ is a linear basis of the
quotient algebra $Lie(X)/Id(S)$.

Some later (see Bokut \cite{b72}, 1972) the Shirshov Composition
lemma was reformulated in the following form: Let $S$ be a  closed
under composition set of monic Lie polynomials  (it means that any
composition $[f,g]_w$ of intersection and inclusion of elements of
$S$ is trivial, i.e., $[f,g]_w=\sum\alpha_i[a_is_ib_i],\
\overline{[a_is_ib_i]}=a_i\overline{s_i}b_i<w,\ a_i,b_i\in X^*,\
s_i\in S,\ \alpha_i\in k$). If $f\in Id(S)$, then $\bar f= a\bar s
b$ for some $s\in S$. And $S$-irreducible NLSW's is a linear basis
of the quotient algebra $Lie(X)/Id(S)$.

The modern form of Shirshov's lemma is the following (see, for
example, Bokut, Chen \cite{BC07}).

{\bf Shirshov's Composition-Diamond Lemma for Lie algebras.} Let
$Lie(X)$ be a free Lie algebra over a field, $S$ monic subset of
$Lie(X)$ relative to some monomial order on $X^*$. Then the
following conditions are equivalent:

1) $S$ is a Gr\"{o}bner-Shirshov basis (i.e., any composition of
intersection and inclusion of elements of $S$ is trivial).

2) If $f\in Id(S)$, then  $\bar f=a\bar sb$ for some $s\in S,\
a,b\in X^*$.

3) $Irr(S)=\{[u]|  [u] \mbox{ is an NLSW and } u \mbox{ does not
contain any } \bar s , s\in S\}$ is a linear basis of the Lie
algebra $Lie(X|S)$ with defining relations $S$.

The proof of the Shirshov's Composition-Diamond lemma for Lie
algebras becomes a model for proofs of number of Composition-Diamond
lemmas for many classes of algebras. An idea of his proof is to
rewrite any element of Lie ideal, generated by $S$ in a form
$$
\sum\alpha_i[a_is_ib_i],
$$
where each $s_i\in S,\ a_i, b_i\in X^*,\ \alpha_i\in k$ such that

i) leading words of each $[a_is_ib_i]$ is equal to $a_i\overline{
s_i}b_i$ (in this case an expression $[asb]$ is called normal Lie
$S$-word in $X$) and

ii) $a_1\overline {s_1}b_1 > a_2\overline {s_2}b_2>\dots.$

Now let $F(X)$ be a free algebra of a variety (or category) of
algebras. Following the idea of  Shirshov's proof, one needs

1) to define appropriate linear basis (normal words) of F(X),

2) to define monomial order of normal words,

3) to define compositions of element of $S$ (they may be
compositions of intersection, inclusion and left (right)
multiplication, or may be else),

4) prove  two key lemmas:

Key Lemma 1. Let $S$ be a Gr\"{o}bner-Shirshov basis (any
composition of polynomials from $S$ is trivial). Then  any $S$-word
is a linear combination of normal $S$-words.

Key Lemma 2. Let $S$ be a Gr\"{o}bner-Shirshov basis, $[a_1s_1b_1]$
and $[a_2s_2b_2]$ normal $S$-words, $s_1, s_2\in S$. If
$a_1\overline{s_1}b_1 = a_2\overline{s_2}b_2$, then $[a_1s_1b_1] -
[a_2s_2b_2]$ is going to zero by elimination of leading words of $S$
(elimination means composition of inclusion).

There are number of CD-lemmas  that realized Shirshov's approach to
them.

Shirshov \cite{Sh62b}, 1962,  assumed implicitly that his approach,
based on the definition of composition of any (not necessary Lie)
polynomials, is equally valid for associative algebras as well (the
first author is a witness that Shirshov understood it very clearly
and explicitly; only lack of non-trivial applications prevents him
from publication this approach for associative algebras). Explicitly
it was done by Bokut \cite{b76} and Bergman \cite{b}.

CD-lemma for associative algebras is formulated and proved in the
same way as for Lie algebras.

{\bf  Composition-Diamond Lemma for associative algebras.} Let
$k\langle X\rangle$ be a free associative algebra  over a field $k$
and a set $X$. Let us fix some monomial order on $X^*$.  Then the
following conditions are equivalent for any monic subset $S$ of
$k\langle X\rangle$:

1) $S$ is a Gr\"{o}bner-Shirshov basis (that is any composition of
intersection and inclusion is trivial).

2) If $f\in Id(S)$, then $\bar f=a\bar sb$ for some $s\in S,\ a,b\in
X^*$.

3) $Irr(S) = \{ u \in X^* |  u \neq a\bar{s}b,s\in S,a ,b \in X^*\}$
is a linear basis of the factor algebra $k\langle
X|S\rangle=k\langle X\rangle/Id(X)$.

\ \

There are a lot of applications of Shirshov's CD-lemmas for Lie and
associative algebras. Let us mention some connected to the Malcev
embedding problem for semigroup algebras (Bokut \cite{bo691,bo692},
1969, there is a semigroup $S$ such that the multiplication
semigroup of the semigroup algebra $k(S)$, where $k$ is a field, is
embeddable into a group, but $k(S)$ is not embeddable into any
division algebra), the unsolvability of the word problem for Lie
algebras (Bokut \cite{b72}), Gr\"{o}bner-Shirshov bases for
semisimple Lie algebras (Bokut, Klein
\cite{BoKl961,BoKl-ICAC97,BoKlJPAA98, BoKl962}), Kac-Moody algebras
(Poroshenko  \cite{Por2000,Por02,Por2002}), finite Coxeter groups
(Bokut, Shiao \cite{BoSh01}), braid groups in different set of
generators (Bokut, Chainikov, Shum \cite{BoChaSh07}, Bokut
\cite{Bo08}, Bokut \cite{Bo09}), quantum algebra of type $A_n$
(Bokut, Malcolmson \cite{BoMa99}), Chinese monoids (Chen, Qiu
\cite{CQ}).

There are applications of Shirshov's CD-lemma \cite{Sh62a}, 1962 for
free anti-commutative non-associative algebras: there are two
anti-commutative Gr\"{o}bner-Shirshov bases of a free Lie algebra,
one gives the Hall basis (Bokut, Chen, Li \cite{BoChLi09}), another
the Lyndon-Shirshov basis (Bokut, Chen, Li \cite{BoChLi10}).

Bokut, Chen, Mo \cite{BoChMo10} proved and reproved some embedding
theorems for associative algebras, Lie algebras, groups, semigroups,
differential algebras, using Shirshov's CD-lemmas for associative
and Lie algebras.

Bahturin, Olshanskii \cite{BahOl} found embeddings without
distortion of associative algebras and Lie algebras into 2-generated
simple algebras. They also used Shirshov's CD-lemmas for associative
and Lie algebras.

Mikhalev \cite{MikhPetr} used Shirshov's approach and CD-lemma for
associative algebras to prove CD-lemma for colored Lie
super-algebras.

Mikhalev, Zolotykh  \cite{MZ} proved CD-lemma for free associative
algebra over a commutative algebra.A Free object in this category is
$k[Y]\otimes k\langle X\rangle$, tensor product of a polynomial
algebra and a free associative algebra. Here one needs to use
several compositions of intersection and inclusion.

Bokut, Fong, Ke \cite{bfk} proved CD-lemma for free associative
conformal (in a sense of V. Kac \cite{Kac96} algebra $C(X, (n),
n=0,1,\dots,D, N(a,b), a,b\in X )$ of a fixed locality $N(a,b)$. A
linear basis of free associative conformal algebra was constructed
by M. Roitman \cite{Ro99}. Any normal conformal word has a form $[u]
= a_1(n_1)[a_2(n_2)[\dots[a_k(n_k)D^ia_{k+1}]\dots]]$, where $a_j\in
X, n_j< N(a_j, a_{j+1}),\ i\geq 0$. The same word without brackets
is called the leading associative word of $[u]$. One needs to use
external multi-operator semigroup as a set of leading associative
words of conformal polynomials (it is the same as for Lie algebras),
several compositions of inclusion and intersection, and new
compositions of left (right) multiplication (last compositions are
absent into classical cases). Also in the CD-lemma  for conformal
algebras we have $1)\Rightarrow2) \Leftrightarrow 3)$, but in
general from 2) does not follow 1). Here conditions 2) and 3) are
formulated in terms of associative leading words, the same as for
Lie algebras. We see that CD-lemma for associative conformal
algebras has a lot of in common with CD-lemma for Lie algebras, but
there are also some differences. Though PBW-theorem is not valid for
Lie conformal algebras (M. Roitman \cite{Ro99}),  some generic
intersection compositions for universal enveloping algebra $U(L)$ of
any Lie conformal algebra $L$  are trivial (it is called ``1/2 PBW
theorem").

Bokut, Chen, Zhang \cite{BCZ} proved CD-lemma for associative
$n$-conformal algebras, where instead of one derivation $D$ and
polynomial algebra $k[D]$ one has $n$ derivations $D_1,\dots,D_n$
and polynomial algebra $k[D_1,\dots,D_n]$. This case is treated in
the same way as for $n=1$. A more general case, the associative
$H$-conformal algebra (or $H$-pseudo-algebra in a sense of Bakalov,
D'Andrea, Kac \cite{BaAnKac01}), where H is any Hopf algebra, is
still open.

Mikhalev, Vasilieva \cite{MV} proved CD-lemma for free
supercommutative polynomial algebras. Here they use compositions of
multiplication as well.

Bokut, Chen, Li \cite{bcli08} proved CD-lemma for free pre-Lie
algebras (also known as Vinberg-Koszul-Gerstenhaber right-symmetric
algebras).

Bokut, Chen, Liu \cite{bcl} proved CD-lemma for free dialgebras in a
sense of Loday \cite{Loday}. Here conditions 1) and 2) are not
equivalent but from 1) follows 2).

The cases of associative conformal algebras and dialgebras show that
definition  of Gr\"{o}bner-Shirshov bases by condition 1) is in
general preferable than the one using 2).

Bokut, Shum \cite{BoShum08} proved CD-lemma for free
$\Gamma$-associative algebras, where $\Gamma$ is a group. It has
applications to the Malcev problem above and to Bruhat normal forms
for algebraic groups.

Eisenbud, Peeva, Sturmfels \cite{EiPeStur} found non-commutative
Gr\"{o}bner basis of any commutative algebra (extending any
commutative Gr\"{o}bner basis to a non-commutative one).

Bokut, Chen, Chen \cite{BCC10} proved CD-lemma for Lie algebras over
commutative algebras. Here one needs to establish Key Lemma 1 in a
more strong form -- any Lie $S$-word is a linear combination of
$S$-words of the form $[asb]_{\bar s}$ in the sense of Shirshov's
special Lie bracketing. As an application they proved Cohn's
conjecture \cite{Conh63} for the case of characteristics 2, 3 and 5
(that some Cohn's examples of Lie algebras over commutative algebras
are not embeddable into associative algebras over the same
commutative algebras).

Bokut, Chen, Deng \cite{bcd} proved CD-lemma for free associative
Rota-Baxter algebras.  As an application, Chen and Mo \cite{CM10}
proved that any dendriform algebra is embeddable into universal
enveloping Rota-Baxter algebra. It was Li Guo's conjecture, \cite
{GL}.

Bokut, Chen, Chen \cite{BCC08} proved CD-lemma for tensor product of
two free associative algebras.  As an application they extended any
Mikhalev-Zolotyh commutative-non-commutative Gr\"{o}bner-Shirshov
basis laying into tensor product $k[Y]\otimes k\langle X\rangle$ to
non-commutative-non-commutative Gr\"{o}bner-Shirshov basis laying
into $k\langle Y\rangle\otimes k\langle X\rangle$ (a la
Eisenbud-Peeva-Sturfels above). They also gave another proof of the
Eisenbud-Peeva-Sturmfekls theorem above.

As we mentioned in introduction, Farkas, Feustel, Green \cite{ffg93}
proved CD-lemma for path algebras.

Drensky, Holtkamp \cite{DH 2008} proved CD-lemma for nonassociative
algebras with multiple linear operators.

Bokut, Chen, Qiu \cite{BCQ} proved CD-lemma for associative algebras
with multiple linear operators.

Dotsenko, Khoroshkin \cite{DK10} proved CD-lemma for operads.

\section{Free categories and free category partial algebras}
Let $X=(V(X),E(X))$ be an oriented (multi)  graph. Then the free
category on $X$ is $C(X)=(Ob(X), Arr(X))$, where $Ob(X)=V(X),$ and
$Arr(X)$ is the set of all paths (``words") of $X$ including the
empty paths $1_{v}$, $v\in V(X)$. It is easy to check $C(X)$ has the
following universal property. Let $\mathscr{C}$ be a category and
$\Gamma_{\mathscr{C}}$ the graph relative to $\mathscr{C}$ i.e.,
$V(\Gamma_{\mathscr{C}})=Ob(\mathscr{C})$ and
$E(\Gamma_{\mathscr{C}})=mor(\mathscr{C})$. Let $e:X\rightarrow
C(X)$ be a  mono graph morphism of the graph $X$ to the graph
$\Gamma_{C(X)},$ where $e=(e_1,e_2)$, and $e_1$ is a mapping on
$V(X),$ $e_2$ on $E(X),$ both $e_1$ and $ e_2$ are mono. For any
graph morphism $b$ from $X$ to $\Gamma_{\mathscr{C}},$ where
$b=(b_1,b_2)$, and $b_1$ is mono, there exists a unique category
morphism (a functor) $f:C(X)\rightarrow \mathscr{C}$, such that the
corresponding diagram is commutative i.e., $fe=b$. Therefore each
category $\mathscr{C}$ is a homomorphic image of a free category
$C(X)$ for some graph $X$ and thus $\mathscr{C}$ is isomorphic to
$C(X)/\rho(S)$ for some set $S$, where $S=\{(u,v)|u,v\mbox{ have\
the\ same\ sources\ and\ the\ same\ targets}\}\subseteq Arr(X)\times
Arr(X)$ and $\rho(S)$ the congruence of $C(X)$ generated by $S$. If
this is the case, $X$ is called the generating set of $\mathscr{C}$
and $S$ the relation set of $\mathscr{C}$ and we denote
$\mathscr{C}=C(X|S)$.

Let $\mathscr{C}$ be a category and $k$ a field. Let
\begin{eqnarray*}
k\mathscr{C}&=&\{f=\sum_{i=1}^n \alpha_i\mu_i|\alpha_i\in k, \
\mu_i\in mor(\mathscr{C}),\ n\geq0, \\
&&\ \  \mu_{i}\ (0\leq i\leq n) \mbox{ have the same domains and the
same codomains} \}.
\end{eqnarray*}

Note that in $k\mathscr{C}$, for $f,\ g\in k\mathscr{C},\ f+g$ is
defined only if $f,\ g$ have the same domain and the same codomain.

A multiplication $\bullet $ in $k\mathscr{C}$ is defined by linearly
extending the usual compositions of morphisms of the category
$\mathscr{C}$. Then ($k\mathscr{C},\ \bullet$) is called the
category partial algebra over $k$ relative to $\mathscr{C}$ and
$kC(X)$ the free category partial algebra generated by the graph
$X$.

\section{Composition-Diamond lemma for categories }
Let $X$ be a oriented (multi) graph, $C(X)$ the free category
generated by $X$ and $kC(X)$ the free category partial algebra.
Since we only consider the morphisms of the free category $C(X)$, we
write $C(X)$ just for $Arr(X)$.

Note that for $f,g\in kC(X)$ if we write $gf$, it means $gf$ is
defined.

A well ordering $>$ on $C(X)$ is called monomial if it satisfies the
following conditions: $u>v\Rightarrow uw>vw$ and $wu>wv$, for any
$u,v,w\in C(X)$. In fact, there are many monomial orders on $C(X)$.
For example, let $E(X)$ be a well ordered set. Then the deg-lex
order $>$ on $C(X)$ is defined by the following way: for any words
$u=x_1\cdots x_m$, $v=y_1\cdots y_n\in$ $C(X),\ m=|u|,\ n=|v|$,
\begin{eqnarray*}
&&u>v\Longleftrightarrow  |u|>|v| \ \mbox { or }\\
&&\ \ \ \ \ \ \ (|u|=|v| \ \ \mbox {and} \ \
x_1=y_1,x_2=y_2,\dots,x_t=y_t,\ x_{t+1}>y_{t+1}\ \ \mbox {for some
}\ 0\leq t<n).
\end{eqnarray*}
It is easy to check that $>$ is a monomial order on $C(X)$. In the
following sections, we will see other monomial orders. Now, we
suppose that $>$ is a fixed monomial order on $C(X)$. Given a
nonzero polynomial $f\in kC(X)$, it has a word $\bar f \in C(X)$
such that $ f=\alpha\overline{f}+\sum{\alpha}_iu_i, $ where
$\overline{f}>u_i, \ 0\neq\alpha , {\alpha}_i \in k, \ u_i\in C(X)$.
We call $\overline{f}$ the leading term of $f$ and $f$ is monic if
$\alpha=1$.

\ \

Let $S\subset kC(X)$ be a set of monic polynomials, $s\in S$ and $u\in C(X)$. We define $S$-word $%
u_s$ by induction:
\begin{enumerate}
\item[(i)] $u_s=s$ is an $S$-word of $s$-length
1.
\item[(ii)] Suppose that $u_s$ is an $S$-word of $s$-length $m$ and $v$
is a word of length $n$, i.e., the number of edges in $v$ is $n$.
Then
$u_sv$\ and\ $vu_s$
are $S$-words of $s$ length $m+n$.
\end{enumerate}

Note that for any $S$-word $u_s=asb$, where $a,b\in C(X)$, we have
$\overline{asb}=a\bar{s}b$.

Let $f,g$ be monic polynomials in $kC(X)$. Suppose that there exist
$w,a,b\in C(X)$ such that $w=\bar f =a\bar g b$. Then we define the
composition of inclusion
\begin{equation*}
(f,g)_{w}=f-agb.
\end{equation*}
For the case that $w=\bar{f}b=a\bar{g}$, $w,a,b\in C(X)$, the
composition of intersection is defined as follows:
\begin{equation*}
(f,g)_{w}=fb-ag.
\end{equation*}

It is clear that
\begin{equation*}
(f,g)_{w}\in Id(f,g)\ \ and\ \ \overline{(f,g)_{w}}<w,
\end{equation*}
where $Id(f,g)$ is the ideal of $kC(X)$ generated by $f,g$.

The composition $(f,g)_w$ is trivial modulo $(S,w)$, if
\begin{equation*}
(f,g)_{w}=\sum\limits_i\alpha_ia_is_ib_i
\end{equation*}
where each $\alpha_i\in k,\ a_i,b_i\in C(X),\ s_i\in S,\ a_is_ib_i$
an $S$-word and $a_i\bar{s_i}b_i<\bar f$. If this is the case, then
we write $(f,g)_{w}\equiv 0\ mod (S,w)$. In general, for $p,q\in
kC(X)$, we write
$$
p\equiv q\quad mod(S,w)
$$
which means that $p-q=\sum\alpha_i a_i s_i b_i $, where each
$\alpha_i\in k,a_i,b_i\in C(X),\ s_i\in S,\ a_is_ib_i$ an $S$-word
and $a_i \bar {s_i} b_i<w$.

\begin{definition}
Let $S\subset kC(X)$ be a nonempty set of monic polynomials. Then
$S$ is called a Gr\"{o}bner-Shirshov basis in $kC(X)$ if any
composition $(f,g)_{w}$ with $f,g\in S$ is trivial modulo $(S,w)$,
i.e., $(f,g)_{w}\equiv0$ $mod(S,w)$.
\end{definition}

\begin{lemma}\label{L1}
Let $a_1s_1b_1,~a_2s_2b_2$ be monic $S$-words. If $S$ is a
Gr\"{o}bner-Shirshov basis in $kC(X)$ and
$w=a_1\overline{s_1}b_1=a_2\overline{s_2}b_2$, then
\begin{equation*}
a_1s_1b_1\equiv a_2s_2b_2\ mod (S,w).
\end{equation*}

\end{lemma}

\noindent\textbf{Proof.} There are three cases to consider.

Case 1. Suppose that subwords $\bar s_1$ and $\bar s_2$ of $w$ are
disjoint, say, $|a_2|\geq |a_1|+|\bar s_1|$. Then, we can assume
that
$
a_2=a_1\bar s_1 c \ \ and  \ b_1=c\bar s_2 b_2
$
for some $c\in C(X)$, and so, $ w=a_1\bar s_1 c \bar s_2 b_2. $ Now,
\begin{eqnarray*}
a_1 s_1 b_1-a_2 s_2 b_2&=&a_1 s_1 c \bar s_2
b_2-a_1\bar s_1 c  s_2 b_2\\
&=&a_1 s_1 c (\bar s_2 - s_2) b_2+a_1(s_1-\bar s_1) c  s_2 b_2.
\end{eqnarray*}
Since $\overline{\overline{s_2}-s_2}<\bar s_2$ and
$\overline{s_1-\overline{s_1}}<\bar s_1$, we conclude that
$$
a_1 s_1 b_1-a_2 s_2 b_2=\sum\limits_i
\alpha_iu_is_1v_i+\sum\limits_j \beta_ju_js_2v_j
$$
for some $\alpha_i,\beta_j\in k$, $S$-words $u_is_1v_i$ and
$u_js_2v_j$ such that $ u_i\bar s_1v_i,u_j\bar s_2v_j<w. $

This shows that $a_1 s_1 b_1\equiv a_2 s_2 b_2 \  \ \ mod(S,w)$.

Case 2. Suppose that the subword $\bar s_1$ of $w$ contains $\bar
s_2$ as a subword. We may assume that $ \bar s_1=a\bar s_2b, \
a_2=a_1a \mbox{ and } b_2=bb_1, \mbox{ that is, } w=a_1a\bar s_2bb_1
$ for some $S$-word $a s_2 b$. We have
\begin{eqnarray*}
a_1 s_1 b_1-a_2 s_2 b_2 &=& a_1 s_1 b_1-a_1 a s_2 b b_1\\
&=& a_1s_1-as_2bb_1\\
&=& a_1(s_1,s_2)_{\overline{s_1}} b_1\\
 &\equiv&0 \  \ \ mod(S,w)
\end{eqnarray*}
since $S$ is a Gr\"{o}bner-Shirshov basis.

Case 3. $\bar{s}_1$ and $\bar{s}_2$ have a nonempty intersection as
a subword of $w$. We may assume that $ a_2=a_1a, \ b_1=bb_2,
w_1=\bar{s}_1b=a\bar{s}_2. $ Then, we have
\begin{eqnarray*}
a_1s_1b_1-a_2s_2b_2
&=&a_1s_1bb_2- a_1as_2b_2\\
&=&a_1(s_1b-as_2)b_2\\
&=&a_1(s_1,s_2)_{w_1}b_2\\
 &\equiv&0 \  \ \ mod(S,w)
\end{eqnarray*}
This completes the proof.\ \ \ \ $\square$

\begin{lemma}\label{L2}
Let $S\subset kC(X)$ be a subset of monic polynomials and
$Irr(S)=\{u\in C(X) |u\ne a\bar s b,\ a,b\in C(X),\ s\in S \}$. Then
for any $f\in kC(X)$,
\begin{equation*}
f=\sum\limits_{u_i\leq \bar f }\alpha_iu_i+
\sum\limits_{a_j\overline{s_j}b_j\leq\bar f}\beta_ja_js_jb_j
\end{equation*}
where each $\alpha_i,\beta_j\in k, \ u_i\in Irr(S)$ and $a_js_jb_j$
an $S$-word.
\end{lemma}
{\bf Proof.} Let $f=\sum\limits_{i}\alpha_{i}u_{i}\in{kC(X)}$, where
$0\neq{\alpha_{i}\in{k}}$ and $u_{1}>u_{2}>\cdots$. If
$u_1\in{Irr(S)}$, then let $f_{1}=f-\alpha_{1}u_1$. If
$u_1\not\in{Irr(S)}$, then there exist some $s\in{S}$ and
$a_1,b_1\in{C(X)}$, such that $\bar f=u_1=a_1\bar{s_1}b_1$. Let
$f_1=f-\alpha_1a_1s_1b_1$. In both cases, we have
$\bar{f_1}<\bar{f}$. Then the result follows from the induction on
$\bar{f}$. \ \ \ \ $\square$

\ \

\begin{theorem}\label{lt1}(Composition-Diamond lemma for
categories)\ Let $S\subset kC(X)$ be a nonempty set of monic
polynomials and $<$ a monomial order on $C(X)$. Let $Id(S)$ be the
ideal of $kC(X)$ generated by $S$. Then the following statements are
equivalent:
\begin{enumerate}
\item [(i)] $S$ is a Gr\"{o}bner-Shirshov basis in $kC(X)$.

\item [(ii)] $f\in Id(S)\Rightarrow \bar f =a\bar s b$ for some $s\in
S$ and $a,b\in C(X)$.

\item [$(ii)'$] $f\in Id(S)\Rightarrow  f = \alpha_1a_1s_1
b_1+\alpha_2a_2s_2b_2+\cdots$, where each $\alpha_i\in k,\ a_is_i
b_i$ is an $S$-word and $a_1\bar s_1b_1>a_2\bar s_2b_2>\cdots$.

\item [(iii)] $Irr(S)=\{u\in C(X) |u\ne a\bar s b\ a,b\in C(X),\ s\in S \}$ is a linear basis of the
partial algebra $kC(X)/Id(S)=kC(X|S)$.
\end{enumerate}
\end{theorem}

{\bf Proof.} $(i)\Rightarrow (ii)$. \ Let $S$ be a
Gr\"{o}bner-Shirshov basis and $0\neq f\in Id(S)$. Then, we have
$$
f=\sum_{i=1}^n\alpha_ia_is_ib_i,
$$
where each $\alpha_i\in k, \ a_i,b_i\in {C(X)}, \ s_i\in S$ and $
a_is_ib_i$ an $S$-word. Let
$$
w_i=a_i\overline{s_i}b_i, \ w_1=w_2=\cdots=w_l>w_{l+1}\geq\cdots,\
l\geq1.
$$
We will use the induction on $l$ and $w_1$ to prove that
$\overline{f}=a\overline{s}b$ for some $s\in S \ \mbox{and} \ a,b\in
{C(X)}$.

If $l=1$, then
$\overline{f}=\overline{a_1s_1b_1}=a_1\overline{s_1}b_1$ and hence
the result holds. Assume that $l\geq 2$. Then, by Lemma \ref{L1}, we
have
$$
a_1s_1b_1\equiv a_2s_2b_2 \ \ mod(S,w_1).
$$
Thus, if $\alpha_1+\alpha_2\neq 0$ or $l>2$, then the result holds
by induction on $l$. For the case $\alpha_1+\alpha_2= 0$ and $l=2$,
we use the induction
on $w_1$. Now, the result follows.\\

$(ii)\Rightarrow (ii)'$. \ Assume (ii) and $0\neq f\in Id(S)$. Let
$f=\alpha_1\overline{f}+\cdots$. Then, by (ii),
$\overline{f}=a_1\overline{s_1}b_1$. Therefore,
$$
f_1=f-\alpha_1a_1s_1b_1, \ \overline{f_1}<\overline{f}, \ f_1\in
Id(S).
$$
Now, by using induction on $\overline{f}$, we have $(ii)'$.\\

$(ii)'\Rightarrow (ii)$. This part is clear.\\

$(ii)\Rightarrow(iii)$. Suppose that $\sum\limits_{i}\alpha_iu_i=0$
in $kC(X|S)$, where $\alpha_i\in k$, $u_i\in {Irr(S)}$. It means
that $\sum\limits_{i}\alpha_iu_i\in{Id(S)}$ in ${kC(X)}$. Then all
$\alpha_i$ must be equal to zero. Otherwise,
$\overline{\sum\limits_{i}\alpha_iu_i}=u_j\in{Irr(S)}$ for some $j$
which contradicts (ii).

Now, by Lemma \ref{L2}, (iii) follows.\\

$(iii)\Rightarrow(i)$. For any $f,g\in{S}$ , by Lemma \ref{L2} and
(iii), we have $ (f,g)_{w}\equiv0\ \ \ mod(S,w). $ Therefore, $S$ is
a Gr\"{o}bner-Shirshov basis. \ \ \ \ $\square$

\ \

\noindent {\bf Remark.} If the category in Theorem \ref{lt1} has
only one object, then Theorem \ref{lt1} is exact Composition-Diamond
lemma for free associative algebras.

\section{Gr\"{o}bner-Shirshov bases for the simplicial category and the cyclic category}

In this section, we give Gr\"{o}bner-Shirshov bases for  the
simplicial category and the cyclic category respectively.

For each non-negative integer $p$, let $[p]$ denote the set
$\{0,1,2,\ldots,p\}$ of integers in their usual order. A (weakly)
monotonic map $\mu:[q]\rightarrow [p]$ is a function on $[q]$ to
$[p]$ such that $i\leq j$ implies $\mu (i)\leq \mu (j)$. The objects
$[p]$ with morphisms all weakly monotonic maps $\mu$ constitute a
category $\mathscr{L}$ called simplicial category. It is convenient
to use two special families of monotonic maps
$$\varepsilon_q^i:[q-1]\rightarrow [q],\ \ \
\eta_q^i:[q+1]\rightarrow [q]$$ defined for $i=0,1,...q$ (and for
$q>0$ in the case of $\varepsilon^i$) by
\begin{equation*}
\varepsilon_q^i(j)=\left\{
\begin{array}{r@{\quad}l}
j,\text{ \ \ \ \   }  & \ \ \ \mbox{ if } \  i>j, \\
j+1, & \ \ \ \mbox{ if } \   i\leq j, \\
\end{array}%
\right.
\end{equation*}

\begin{equation*}
\eta_q^i(j)=\left\{
\begin{array}{r@{\quad}l}
j,\text{ \ \ \ \   }  & \ \ \ \mbox{ if } \  i\geq j, \\
j-1, & \ \ \ \mbox{ if } \   j>i. \\
\end{array}%
\right.
\end{equation*}

Let $X=(V(X),E(X))$ be an oriented (multi)  graph, where
$V(X)=\{[p]\ |\ p\in Z^+\cup\{0 \}\}$ and $E(X)=\{
\varepsilon_p^i:[p-1]\rightarrow [p], \ \eta_q^j:[q+1]\rightarrow
[q]\ |\ p>0, 0\leq i\leq p, 0\leq j \leq q\}. $ Let $S\subseteq
C(X)\times C(X)$ be the relation set consisting of the following:
\begin{eqnarray*}
&&f_{_{q+1,q}}:\ \ \
\varepsilon_{q+1}^i\varepsilon_{q}^{j-1}=\varepsilon_{q+1}^j\varepsilon_{q}^{i},\
\ \ j>i,\\
&&g_{_{q,q+1}}:\ \ \
\eta_{q}^{j}\eta_{q+1}^{i}=\eta_{q}^i\eta_{q+1}^{j+1},\ \ \
j\geq i,\\
&& h_{_{q-1,q}}:\ \ \ \eta_{q-1}^{j}\varepsilon_{q}^i=\left\{
\begin{array}{r@{\quad}l}
\varepsilon_{q-1}^i\eta_{q-2}^{j-1}, & \ \  \  j>i, \\
1_{q-1}, & \ \  \   i=j, \ \ i=j+1, \\
\varepsilon_{q-1}^{i-1}\eta_{q-2}^{j}, & \ \  \   i>j+1.%
\end{array}%
\right.
\end{eqnarray*}

Then the simplicial category $\mathscr{L}$ is just the category
$C(X|S)$ generated by $X$ with defining relation $S$, see Maclane
\cite{Sa}, Theorem VIII. 5.2. We will give another proof in what
follows.

We order $C(X)$ by the following way.

Firstly, for any $\eta_{p}^{i},\eta_{q}^{j}\in$
$\{\eta_{p}^{i}|p\geq 0, 0\leq i\leq p\}$,
$\eta_{p}^{i}>\eta_{q}^{j}$ iff $p>q$ or ($p=q $ and $i<j$).

Secondly, for each $u=\eta_{p_1}^{i_1}\eta_{p_2}^{i_2}\cdots
\eta_{p_n}^{i_n}\in $ $\{\eta_{p}^{i}|p\geq 0, 0\leq i\leq p\}$*
(all possible words on $\{\eta_{p}^{i}|p\geq 0, 0\leq i\leq p\}$,
including the empty word $1_v$, $v\in Ob(X)$), let
$wt(u)=(n,\eta_{p_n}^{i_n},\eta_{p_{n-1}}^{i_{n-1}},\cdots,\eta_{p_1}^{i_1}).$
Then for any $u,v\in$ $\{\eta_{p}^{i}|p\geq 0, 0\leq i\leq p\}$*,
$u>v$ iff $wt(u)>wt(v)$ lexicographically.

Thirdly, for any $\varepsilon_{p}^{i},\varepsilon_{q}^{j}\in$
$\{\varepsilon_p^{i},|p\in Z^+, 0\leq i\leq p\}$,
$\varepsilon_{p}^{i}>\varepsilon_{q}^{j}$ iff $p>q$ or ($p=q $ and
$i<j$).

Finally, for each
$u=v_{0}\varepsilon_{p_{1}}^{i_1}v_1\varepsilon_{p_{2}}^{i_2}\cdots\varepsilon_{p_{n}}^{i_n}
v_{n}\in C(X),$ $n\geq 0$, $v_j\in \{\eta_{p}^{i}|p\geq 0, 0\leq
i\leq p\}$*, let
$wt(u)=(n,v_0,v_1,\cdots,v_n,\varepsilon_{p_{1}}^{i_1},\cdots,\varepsilon_{p_{n}}^{i_n}).$
Then for any $u,v\in C(X)$,
$$
 u\succ_{_1}v \Leftrightarrow wt(u)>wt(v)\ \ \ \mbox{lexicographically}.
$$

It is easy to check that the $\succ_{_1}$ is a monomial order on
$C(X)$. Then we have the following theorem.

\begin{theorem}
Let $X$, $S$ be defined as the above, the generating set and the
relation set of the quotient category $C(X|S)$ respectively. Then
with the order $\succ_1$ on $C(X)$, $S$ is  a Gr\"{o}bner-Shirshov
basis for the category partial algebra $kC(X|S)$.
\end{theorem}

{\bf Proof.} According to the order $\succ_{_1}$,
${\bar{f}_{_{q+1,q}}}=\varepsilon_{q+1}^i\varepsilon_{q}^{j-1}$,
${\bar{g}_{_{q,q+1}}}=\eta_{q}^{j}\eta_{q+1}^{i}$ and
${\bar{h}_{_{q-1,q}}}=\eta_{q-1}^{j}\varepsilon_{q}^i$. So, all the
possible compositions of $S$ are the following:
\begin{enumerate}
\item[(a)]$(f_{_{q+2,q+1}},f_{_{q+1,q}})_{_{\varepsilon_{q+2}^k\varepsilon_{q+1}^i\varepsilon_{q}^{j-1}}}$,
\ \ $k\leq i\leq j-1$;
\item[ (b)]
$(g_{_{q-1,q}},g_{_{q,q+1}})_{_{\eta_{q-1}^{k}\eta_{q}^{j}\eta_{q+1}^{i}}}$,
\ \ $i\leq j\leq k$;
\item[ (c)]
$(h_{_{q,q+1}},f_{_{q+1,q}})_{_{\eta_{q}^{k}\varepsilon_{q+1}^i\varepsilon_{q}^{j-1}}}$,
\ \ $i\leq j-1$;
\item[(d)]
$(g_{_{q-2,q-1}},h_{_{q-1,q}})_{_{\eta_{q-2}^{k}\eta_{q-1}^{j}\varepsilon_{q}^{i}}}$,
\ \ $j\leq k$.
\end{enumerate}

We will prove that all possible compositions are trivial. Here, we
only give the proof of the (b). For others cases, the proofs are
similar.

Let us consider the following subcases of the case (b): (I) $i<j<k$;
(II) $i<j,j=k, \mbox{or} \ j=k+1$; (III) $i<k,k+1<j$; (IV)
$j>k+1,i=k, k+1$; (V) $j>i>k+1$.

For subcase (I),
\begin{eqnarray*}
(h_{_{q,q+1}},f_{_{q+1,q}})_{_{\eta_{q}^{k}\varepsilon_{q+1}^i\varepsilon_{q}^{j-1}}}
&=&\varepsilon_q^i\eta_{q-1}^{k-1}\varepsilon_q^{j-1}-\eta_q^
k\varepsilon_{q+1}^j\varepsilon_q^i\\
&\equiv& \varepsilon_q^i\varepsilon_{q-1}^{j-1}\eta_{q-2}^{k-2}-\varepsilon_{q}^j\eta_{q-1}^{k-1}\varepsilon_q^i\\
&\equiv&0 \  \ \
mod(S,{\eta_{q}^{k}\varepsilon_{q+1}^i\varepsilon_{q}^{j-1}}).
\end{eqnarray*}

For subcase (II),
\begin{eqnarray*}
(h_{_{q,q+1}},f_{_{q+1,q}})_{_{\eta_{q}^{k}\varepsilon_{q+1}^i\varepsilon_{q}^{j-1}}}
&=&\varepsilon_q^i\eta_{q-1}^{k-1}\varepsilon_q^{j-1}-\eta_q^k\varepsilon_{q+1}^j\varepsilon_q^i\\
&\equiv&0
 \  \ \
mod(S,{\eta_{q}^{ k}\varepsilon_{q+1}^i\varepsilon_{q}^{j-1}}).
\end{eqnarray*}

For subcase (III),
\begin{eqnarray*}
(h_{_{q,q+1}},f_{_{q+1,q}})_{_{\eta_{q}^{k}\varepsilon_{q+1}^i\varepsilon_{q}^{j-1}}}
&=&\varepsilon_q^i\eta_{q-1}^{k-1}\varepsilon_q^{j-1}-\eta_q^k\varepsilon_{q+1}^j\varepsilon_q^i\\
&\equiv& \varepsilon_q^i\varepsilon_{q-1}^{j-2}\eta_{q-2}^{k-1}-\varepsilon_{q}^{j-1}\eta_{q-1}^{k}\varepsilon_q^i\\
&\equiv& 0 \  \ \
mod(S,{\eta_{q}^{k}\varepsilon_{q+1}^i\varepsilon_{q}^{j-1}}).
\end{eqnarray*}

For subcase (IV),
\begin{eqnarray*}
(h_{_{q,q+1}},f_{_{q+1,q}})_{_{\eta_{q}^{k}\varepsilon_{q+1}^i\varepsilon_{q}^{j-1}}}
&=&\varepsilon_q^{j-1}-\eta_q^k\varepsilon_{q+1}^j\varepsilon_q^i\\
&\equiv& \varepsilon_q^{j-1}-\varepsilon_q^{j-1}\eta_{q-1}^k\varepsilon_q^i\\
&\equiv&0 \  \ \
mod(S,{\eta_{q}^{k}\varepsilon_{q+1}^i\varepsilon_{q}^{j-1}}).
\end{eqnarray*}

For subcase (V),
\begin{eqnarray*}
(h_{_{q,q+1}},f_{_{q+1,q}})_{_{\eta_{q}^{k}\varepsilon_{q+1}^i\varepsilon_{q}^{j-1}}}
&=&\varepsilon_q^{i-1}\eta_{q-1}^{k}\varepsilon_q^{j-1}-\eta_q^k\varepsilon_{q+1}^j\varepsilon_q^i\\
&\equiv& \varepsilon_q^{i-1}\varepsilon_{q-1}^{j-2}\eta_{q-2}^{k}-\varepsilon_{q}^{j-1}\eta_{q-1}^{k}\varepsilon_q^i\\
&\equiv&0 \  \ \
mod(S,{\eta_{q}^{k}\varepsilon_{q+1}^i\varepsilon_{q}^{j-1}}).
\end{eqnarray*}

Therefore $S$ is a Gr\"{o}bner-Shirshov basis of the category
partial algebra $kC(X|S)$.\ \ \ \ $\square$

\ \

By Theorem \ref{lt1}, $
Irr(S)=\{\varepsilon_p^{i_1}...\varepsilon_{p-m+1}^{i_m}\eta_{q-n}^{j_1}...\eta_{q-1}^{j_n}|p\geq
i_1>...> i_m\geq 0, \ \ 0\leq j_1< ... < j_n <q, \mbox{and }
q-n+m=p\} $ is a linear basis of the category partial algebra
$kC(X|S)$. Therefore, we have the following corollaries.

\begin{corollary}(Maclane \cite{Sa}, Lemma VIII. 5.1)
In the category $C(X|S)$, each morphism $\mu:[q]\rightarrow [p]$ can
be uniquely represented as
$$
\varepsilon_p^{i_1}...\varepsilon_{p-m+1}^{i_m}\eta_{q-n}^{j_1}...\eta_{q-1}^{j_n},$$
where $p\geq i_1>...> i_m\geq 0, \ \ 0\leq j_1< ... < j_n <q,
\mbox{and } q-n+m=p$.
\end{corollary}

\begin{corollary}(Maclane \cite{Sa}, Theorem VIII. 5.2)
$\mathscr{L}=C(X|S)$.
\end{corollary}

The cyclic category is defined by generators and defining relations
as follows, see \cite{G}. Let $Y=(V(Y),E(Y))$ be an oriented (multi)
graph, where $V(Y)=\{[p]\ |\ p\in Z^+\cup\{0 \}\}$, and $E(Y)=\{
\varepsilon_p^i:[p-1]\rightarrow [p], \ \eta_q^j:[q+1]\rightarrow
[q],\ t_q:[q]\rightarrow [q] |\ p>0, 0\leq i\leq p, 0\leq j \leq
q\}.$ Let $S\subseteq C(Y)\times C(Y)$ be the set consisting of the
following relations:
\begin{eqnarray*}
&&f_{_{q+1,q}}:\ \ \
\varepsilon_{q+1}^i\varepsilon_{q}^{j-1}=\varepsilon_{q+1}^j\varepsilon_{q}^{i},\
\ \ j>i,\\
&&g_{_{q,q+1}}:\ \ \
\eta_{q}^{j}\eta_{q+1}^{i}=\eta_{q}^i\eta_{q+1}^{j+1},\ \ \
j\geq i,\\
&& h_{_{q-1,q}}:\ \ \ \eta_{q-1}^{j}\varepsilon_{q}^i=\left\{
\begin{array}{r@{\quad}l}
\varepsilon_{q-1}^i\eta_{q-2}^{j-1}, & \ \  \  j>i, \\
1_{q-1}, & \ \  \   i=j, \ \ i=j+1, \\
\varepsilon_{q-1}^{i-1}\eta_{q-2}^{j}, & \ \  \   i>j+1,%
\end{array}%
\right.\\
&&\rho_1:\ \ \ \ \ \
t_q\varepsilon_{q}^{i}=\varepsilon_{q}^{i-1}t_{q-1}, \ \ \
i=1,...,q,
\\
&&\rho_2:\ \ \ \ \ \ t_q\eta_{q}^{i}=\eta_{q}^{i-1}t_{q+1},\ \ \i=1,...,q, \\
&&\rho_3:\ \ \ \ \ \ t_q^{q+1}=1_q.
\end{eqnarray*}

The category  $C(Y|S)$ is called cyclic category, denoted by
$\Lambda$.

An order on $C(Y)$ is defined by the following way.

Firstly, for any $t_p^i$, $t_q^j\in$ $\{t_q|q\geq 0\}^{*}$,
$(t_p)^i> (t_q)^j$ iff $i>j$ or $(i=j \mbox{ and } p>q)$.

Secondly, for any $\eta_{p}^{i},\eta_{q}^{j}\in$
$\{\eta_{p}^{i}|p\geq 0, 0\leq i\leq p\}$,
$\eta_{p}^{i}>\eta_{q}^{j}$ iff $p>q$ or ($p=q $ and $i<j$).

Thirdly, for each $u=w_0\eta_{p_1}^{i_1}w_1\eta_{p_2}^{i_2}\cdots
w_{n-1} \eta_{p_n}^{i_n}w_n\in $ $\{t_q,\eta_{p}^{i}|q,p\geq 0,
0\leq i\leq p\}$*,where $w_i\in\{t_q|q\geq 0\}^*$, let
$wt(u)=(n,w_0,w_1,\cdots,w_n,\eta_{p_n}^{i_n},\eta_{p_{n-1}}^{i_{n-1}},\cdots,\eta_{p_1}^{i_1}).$
Then for any $u,v\in$ $\{t_q,\eta_{p}^{i}|q,p\geq 0, 0\leq i\leq
p\}$*, $u>v$ iff $wt(u)>wt(v)$ lexicographically.

Fourthly, for any $\varepsilon_{p}^{i},\varepsilon_{q}^{j}\in$
$\{\varepsilon_p^{i},|p\in Z^+, 0\leq i\leq p\}$,
$\varepsilon_{p}^{i}>\varepsilon_{q}^{j}$ iff $p>q$ or ($p=q $ and
$i<j$).

Finally, for each
$u=v_{0}\varepsilon_{p_{1}}^{i_1}v_1\varepsilon_{p_{2}}^{i_2}\cdots\varepsilon_{p_{n}}^{i_n}
v_{n}\in C(Y),$ $n\geq 0$, $v_j\in \{t_q,\eta_{p}^{i}|q,p\geq 0,
0\leq i\leq p\}$*, let
$wt(u)=(n,v_0,v_1,\cdots,v_n,\varepsilon_{p_{1}}^{i_1},\cdots,\varepsilon_{p_{n}}^{i_n}).$

Then for any $u,v\in C(Y)$,
$$
 u\succ_{_2}v \Leftrightarrow wt(u)>wt(v)\ \ \ \mbox{lexicographically}.
$$

It is also easy to check the order $\succ_{_2}$  is a monomial order
on $C(Y)$, which is an extension of $\succ_{_1}$. Then we have the
following theorem.

\begin{theorem}\label{t2}
Let $Y$, $S$ be defined as the above, the generating set and the
relation set of cyclic category $C(Y|S)$ respectively. Let
$S^{C}=S\cup \{\rho_4,\rho_5\}$, where
\begin{eqnarray*}
&\rho_4:& \ \ t_q\varepsilon_{q}^0=\varepsilon_{q}^{q},\\
&\rho_5:&\ \ t_q\eta_{q}^0=\eta_{q}^{q}t_{q+1}^{2}.
\end{eqnarray*}
Then
\begin{enumerate}
\item [(1)] With the order   $\succ_{_2}$ on $C(Y)$, $S^{C}$ is a Gr\"{o}bner-Shirshov
basis for the cyclic category partial algebra $kC(Y|S)$.
\item [(2)] For each morphism $\mu:[q]\rightarrow [p]$ in the cyclic category $\Lambda=C(Y|S)$, $\mu$ can be
uniquely represented as
$$
\varepsilon_p^{i_1}...\varepsilon_{p-m+1}^{i_m}\eta_{q-n}^{j_1}...\eta_{q-1}^{j_n}t_{q}^k,$$
where $p\geq i_1>...> i_m\geq 0, \ \ 0\leq j_1< ... < j_n <q, \ 0
\leq k \leq q \ \mbox{and } \ q-n+m=p.$ \end{enumerate}
\end{theorem}

{\bf Proof.} It is easy to check that
${\bar{f}_{_{q+1,q}}}=\varepsilon_{q+1}^i\varepsilon_{q}^{j-1}$,
${\bar{g}_{_{q,q+1}}}=\eta_{q}^{j}\eta_{q+1}^{i}$,
${\bar{h}_{_{q-1,q}}}=\eta_{q-1}^{j}\varepsilon_{q}^i$,
$\bar{\rho_1}=t_q\varepsilon_{q}^i$, $\bar{\rho_2}=t_q\eta_{q}^i$,
$\bar{\rho_3}=t_{q}^{q+1}$, $\bar{\rho_4}=t_q\varepsilon_{q}^0$, and
$\bar{\rho_5}=t_q\eta_{q}^0$.

First of all, we prove $Id(S)=Id(S^C)$. It suffices to show
$\rho_4,\rho_5\in$ $Id(S)$. Since
$(\rho_3,\rho_1)_{t_q^{q+1}\varepsilon_q^q}$=$t_q^q\varepsilon_q^{q-1}t_{q-1}-\varepsilon_q^q\equiv
t_q\varepsilon_q^0t_{q-1}^q-\varepsilon_q^q\equiv
t_q\varepsilon_q^0-\varepsilon_q^q=\rho_4$ and
$(\rho_3,\rho_2)_{t_q^{q+1}\eta_q^q}$=$t_q^q\eta_q^{q-1}t_{q+1}-\eta_q^q\equiv
t_q\eta_q^0t_{q+1}^q-\eta_q^q$, $\rho_4 \ \ \mbox{and}\ \
t_q\eta_q^0t_{q+1}^q-\eta_q^q \in Id(S)$. Clearly, the leading term
of the polynomial $t_q\eta_q^0t_{q+1}^q-\eta_q^q $ is
$t_q\eta_q^0t_{q+1}^q$. Therefore
$(t_q\eta_q^0t_{q+1}^q-\eta_q^q,\rho_3)_{t_q\eta_q^0t_{q+1}^{q+2}}=-t_q\eta_q^0+\eta_q^qt_{q+1}^2$
and thus $\rho_5\in Id(S).$

Secondly, we prove that all possible compositions of $S^{C}$ are
trivial which are the following:
\begin{enumerate}
\item[(a)]$(f_{_{q+2,q+1}},f_{_{q+1,q}})_{_{\varepsilon_{q+2}^k\varepsilon_{q+1}^i\varepsilon_{q}^{j-1}}}$,
\ \ $k\leq i\leq j-1$;
\item[(b)]
$(g_{_{q-1,q}},g_{_{q,q+1}})_{_{\eta_{q-1}^{k}\eta_{q}^{j}\eta_{q+1}^{i}}}$,
\ \ $i\leq j\leq k$;
\item[(c)]
$(h_{_{q,q+1}},f_{_{q+1,q}})_{_{\eta_{q}^{k}\varepsilon_{q+1}^i\varepsilon_{q}^{j-1}}}$,
\ \ $i\leq j-1$;
\item[(d)]
$(g_{_{q-2,q-1}},h_{_{q-1,q}})_{_{\eta_{q-2}^{k}\eta_{q-1}^{j}\varepsilon_{q}^{i}}}$,
\ \ $j\leq k$;
\item[(e)]$(\rho_1,f_{_{q+1,q}})_{_{t_{q+1}\varepsilon
_{q+1}^{i}\varepsilon_{q}^{j-1}}}$, \ \ $j>i $ and $i=1,2,\dots,q$;
\item[(f)]$(\rho_3,\rho_1)_{{t_{q}^{q+1}}\varepsilon_{q}^{i}}$, \ \ $i=1,2,\dots,q$;
\item[(g)]$(\rho_2,g_{q,q+1})_{{t_{q}}\eta_{q}^{j}\eta_{q+1}^{i}}$, \ \ $j\geq i$ and $j=1,2,\dots,q$;
\item[(h)]$(\rho_2,h_{q,q+1})_{{t_{q}}\eta_{q}^{j}\varepsilon_{q+1}^{i}}$, \ \ $j\geq i$ and $j=1,2,\dots,q$;
\item[(i)]$(\rho_3,\rho_2)_{{t_{q}^{q+1}}\eta_{q}^{i}}$, \ \ $i=1,2,\dots,q$;
\item[(j)]$(\rho_3,\rho_4)_{{t_{q}^{q+1}}\varepsilon_{q}^{0}}$;
\item[(k)]$(\rho_3,\rho_5)_{{t_{q}^{q+1}}\eta_{q}^{0}}$;
\item[(l)]$(\rho_4,f_{_{q+1,q}})_{_{t_{q+1}\varepsilon
_{q+1}^{0}\varepsilon_{q}^{j-1}}}$, \ \ $j>0 $;
\item[(m)]
$(\rho_5,g_{_{q,q+1}})_{_{t_{q}\eta_{q}^{0}\eta_{q+1}^{0}}}$;
\item[(n)]$(\rho_5,h_{_{q,q+1}})_{_{t_{q}\eta_{q}^{0}\varepsilon_{q+1}^{i}}}$, $i\geq 0$.
\end{enumerate}

Here, we only give the proof of the case (n)
$(\rho_6,h_{_{q,q+1}})_{_{t_{q}\eta_{q}^{0}\varepsilon_{q+1}^{i}}}$.
The others can be similarly proved. Let us consider the following
subcases of the case (n): (I) $i=0$; (II) $i=1$; (III) $i>1$.

For subcase (I),
\begin{eqnarray*}
(\rho_6,h_{_{q,q+1}})_{_{t_{q}\eta_{q}^{0}\varepsilon_{q+1}^{0}}}
&=&\eta_q^qt_{q+1}^{2}\varepsilon_{q+1}^0-t_q\\
&\equiv& \eta_q^qt_{q+1}\varepsilon_{q+1}^{q+1}-t_q\\
&\equiv& \eta_q^q\varepsilon_{q+1}^{q}t_{q}-t_q\\
&\equiv& 0 \  \ \ mod(S,t_{q}\eta_{q}^{0}\varepsilon_{q+1}^{0}).
\end{eqnarray*}

For subcase (II),
\begin{eqnarray*}
(\rho_6,h_{_{q,q+1}})_{_{t_{q}\eta_{q}^{0}\varepsilon_{q+1}^{1}}}
&=&\eta_q^qt_{q+1}^{2}\varepsilon_{q+1}^1-t_q\\
&\equiv& \eta_q^qt_{q+1}\varepsilon_{q+1}^{0}t_q-t_q\\
&\equiv& \eta_q^q\varepsilon_{q+1}^{q+1}t_{q}-t_q\\
&\equiv&0 \  \ \ mod(S,t_{q}\eta_{q}^{0}\varepsilon_{q+1}^{1}).
\end{eqnarray*}

For subcase (III),
\begin{eqnarray*}
(\rho_6,h_{_{q,q+1}})_{_{t_{q}\eta_{q}^{0}\varepsilon_{q+1}^{i}}}
&=&\eta_q^qt_{q+1}^{2}\varepsilon_{q+1}^i-t_q\varepsilon_{q}^{i-1}\eta_{q-1}^0\\
&\equiv& \eta_q^q\varepsilon_{q+1}^{i-2}t_{q}^2-\varepsilon_{q}^{i-2}t_{q-1}\eta_{q-1}^0\\
&\equiv& 0 \  \ \ mod(S,t_{q}\eta_{q}^{0}\varepsilon_{q+1}^{i}).
\end{eqnarray*}

Thus $S^C$ is a Gr\"{o}bner-Shirshov basis  of the category partial
algebra $kC(Y|S)$.

Now, by Theorem \ref{lt1}, for each morphism $\mu:[q]\rightarrow
[p]$ in $\Lambda=C(Y|S)$ can be uniquely represented as
$$
\varepsilon_p^{i_1}...\varepsilon_{p-m+1}^{i_m}\eta_{q-n}^{j_1}...\eta_{q-1}^{j_n}t_{q}^k,$$
where $p\geq i_1>...> i_m\geq 0, \ \ 0\leq j_1< ... < j_n <q, \ 0
\leq k \leq q \ \mbox{and } \ q-n+m=p.$ \ \ \ \ $\square$

\ \

\noindent {\bf Remark.} According to Loday \cite{Loday-}, the
uniqueness property in Theorem \ref{t2} (2) was known.

\end{document}